\documentclass[11pt, reqno]{amsart}
\RequirePackage{amsmath,amsthm,amssymb,color,tikz}
\RequirePackage{epsfig,amsfonts,latexsym, hyperref, bbm, mathrsfs}
\RequirePackage{natbib}

\newcommand{\surv}{\mathcal{S}}

\newcommand{\PP}{\mathcal{P}}

\newcommand{\h}{\hspace{2cm}}

\newcommand{\summation}{\sum_{i=1}^{\infty}}
\newcommand{\summationl}{\sum_{l=1}^{\infty}}

\newcommand{\E}{\mathbb{E}}

\newcommand{\T}{\mathbb{T}}
\newcommand{\uv}{\mathbf{v}}
\newcommand{\uu}{\mathbf{e}}

\newcommand{\1}{\mathbbm{1}}

\DeclareMathOperator{\estar}{\mathbf{E}^{*}}
\DeclareMathOperator{\pstar}{\mathbf{P}^{*}}

\DeclareMathOperator{\PT}{\mathbf{P}_{\T}^*}
\DeclareMathOperator{\prob}{\mathbf{P}}
\DeclareMathOperator{\exptn}{\mathbf{E}}

\DeclareMathOperator{\bm}{\mathbbm{1}}

\newcommand{\Xu}{X_{\uu}}
\newcommand{\Sv}{S_{\uv}}
\newcommand{\Iv}{I_{\uv}}

\newcommand{\nkb}{\tilde{N}_n^{(K,B)}}

\newcommand{\tnk}{\tilde{N}_n^{(K)}}
\newcommand{\tnkb}{\tilde{N}_n^{(K,B)}}
\newcommand{\dnib}{D_{n-i}^{(B)}}
\newcommand{\aub}{A_\uu^{(B)}}

\numberwithin{equation}{section}

\newtheorem{thm}{Theorem}[section]
\newtheorem{remark}[thm]{Remark}

\newtheorem{example}{Example}[section]
\newtheorem{defn}[thm]{Definition}

\newtheorem{lemma}[thm]{Lemma}

\newcommand{\eql}{\stackrel{\mathcal{L}}{=}}
\newcommand{\vconv}{\stackrel{v}{\longrightarrow}}
\newcommand{\lfb}{\left(}
\newcommand{\rfb}{\right)}
\newcommand{\lsb}{ \left\lbrace }
\newcommand{\rsb}{\right\rbrace }

\newcommand{\inv}{{-1}}
\newcommand{\beq}{\begin{equation}}
\newcommand{\eeq}{\end{equation}}
\newcommand{\alns}[1]{\begin{align*}#1\end{align*}}
\newcommand{\aln}[1]{\begin{align} #1 \end{align}}
\newcommand{\been}{\begin{enumerate}}
\newcommand{\een}{\end{enumerate}}
\newcommand{\dnk}{D_{n-K}}

\newcommand{\bft}{\mathbf{t}}
\newcommand{\cardi}[1]{| #1 |_i}

\newcommand{\tob}{\T_1^{(B)}}
\newcommand{\caltkb}{\mathcal{T}_{K,B}}

\allowdisplaybreaks

\begin{document}

\allowdisplaybreaks

\title[Point Process convergence for Branching random walks]{Point process convergence for branching random walks with regularly varying steps}

\author[A. Bhattacharya]{Ayan Bhattacharya} \thanks{Ayan Bhattacharya's research was partially supported by the project RARE-318984 (a Marie Curie FP7 IRSES Fellowship).}
\address{Statistics and Mathematics Unit, Indian Statistical Institute, 203 B. T. Road, Kolkata 700108.}

\email{ayanbhattacharya.isi@gmail.com}

\author[R. S. Hazra]{Rajat Subhra Hazra}\thanks{Rajat Subhra Hazra's research was supported by Cumulative Professional Development Allowance from Ministry of Human Resource Development, Government of India and Department of Science and Technology, Inspire funds.}
\email{rajatmaths@gmail.com}

\author[P. Roy]{Parthanil Roy}\thanks{Parthanil Roy's research was supported by Cumulative Professional Development Allowance from Ministry of Human Resource Development, Government of India and the project RARE-318984 (a Marie Curie FP7 IRSES Fellowship).}
\email{parthanil.roy@gmail.com}



\begin{abstract}
We consider the limiting behaviour of the point processes associated with a branching random walk with supercritical branching mechanism and balanced regularly varying step size. Assuming that the underlying branching process satisfies Kesten-Stigum condition, it is shown that the point process sequence of properly scaled displacements coming from the $n^{th}$ generation converges weakly to a Cox cluster process. In particular, we establish that a conjecture of \cite{brunet:derrida:2011} remains valid in this setup, investigate various other issues mentioned in their paper and recover the main result of \cite{durrett:1983} in our framework.
\end{abstract}

\maketitle

\section{{\textbf{Introduction}}}

Suppose $\lsb Z_i \rsb_{i \ge 0}$ is a supercritical Galton-Watson process with $Z_0 \equiv 1$ (the root), branching random variable $Z_1$ (size of the first generation), and progeny mean $\mu :=\exptn (Z_1) \in (1,\infty)$. It is well-known that  $ Z_n / \mu^n$ is a martingale sequence that converges almost surely to a non-negative random variable $W$. We assume that the branching random variable satisfies the \emph{Kesten-Stigum condition}
\begin{equation} \label{nmr1}
\exptn ( Z_1 \log^+ Z_1) < \infty.
\end{equation}
We shall condition on the survival of this Galton-Watson process and \eqref{nmr1} ensures that the limiting random variable $W$ is almost surely positive; see \cite{kesten:stigum:1966a}. This rooted infinite Galton-Watson tree will be denoted in this article by $\T = (V,E)$, where the collection of all vertices is denoted by $V$, the collection of all edges is denoted by $E$ and the root is denoted by $o$. Note that every vertex $\uv$ is connected to the root  $o$ by a unique geodesic path which will be denoted by $\Iv$ and the length of this path will be denoted by $|\uv|$.

We define a branching random walk with balanced regularly varying step size as follows. After obtaining the entire infinite tree $\T$, we assign independent and identically distributed random variables $\{\Xu: \uu \in E\}$ (that are also independent of the Galton-Watson process $\{Z_i\}_{i \geq 0}$) on the edges satisfying the \emph{regular variation condition}
\begin{equation} \label{nmr3}
\prob (|\Xu| > x) = x^{-\alpha} L(x),
\end{equation}
where $\alpha > 0$ and $L(x)$ is a slowly varying function (i.e., for all $x>0$, $L(tx)/L(t) \to 1$ as $t \to \infty$), and the \emph{tail balance condition}
\begin{equation} \label{nmr4}
\frac{\prob (\Xu > x)}{\prob (|\Xu|>x)} \to p \h  \mbox{ and } \h \frac{\prob (\Xu < -x)}{\prob (|\Xu| > x)} \to q
\end{equation}
as $x \to \infty $ for some $p,q \ge 0$ with $p+q=1$. For an encyclopaedic treatment of
regularly varying and slowly varying functions, see \cite{bingham:goldie:teugels:1987}. To each vertex $\uv$, we assign displacement labels $\Sv$, which is the sum of all edge random variables on the geodesic path from the root $o$ to the vertex $\uv$, i.e.,
\begin{equation} \label{nmr2}
\Sv = \sum_{\uu \in \Iv} X_\uu.
\end{equation}
The collection of displacement random variables $\lsb \Sv: |\uv|=n \rsb$ forms the $n^{th}$ generation of our branching random walk.

Branching random walk has been of interest starting from the classical works of \cite{hammersley1974}, \cite{kingman1975}, \cite{biggins1976, biggins:1977a, biggins:1977b}. Recently, extremes of branching random walk has gained much prominence  due to its connection to tree indexed random walk and Gaussian free field; see \cite{BZ2012}, \cite{hu:shi:2009, addario2009}, \cite{aidekon2013}, \cite{madaule2011},
\cite{biskup:louidor:2013, biskup:louidor:2014}, \cite{bramson:ding:zeitouni:2013}. See also \cite{bramson:1978, bramson1983}, \cite{lalley:sellke:1987}, \cite{arguin2011genealogy, arguin2012poissonian, arguin2013extremal}, \cite{aidekon2013branching} for related results on extremes of branching Brownian motion.

Heavy tailed edge random variables were introduced in branching random walks by \cite{durrett:1979, durrett:1983}; see also \cite{kyprianou:1999}, \cite{gantert2000}, and the recent works of \cite{lalley:shao:2013} and \cite{berard2014}. It was shown in \cite{durrett:1983} that when the step sizes have regularly varying tails, then the maximum displacement grows exponentially and converges (after scaling) to a $W$-mixture of Fr\'{e}chet random variables. This limiting behaviour is very different from the ones obtained by \cite{biggins1976} and \cite{bramson:1978} in the light tailed case.

It was predicted in \cite{brunet:derrida:2011} that the limits of point processes of properly normalized displacements of branching random walk and branching Brownian motion should be decorated Poisson point processes. This conjecture was proved to be true for branching Brownian motion by \cite{arguin2012poissonian, arguin2013extremal} and \cite{aidekon2013branching}, and for branching random walks with step sizes having finite exponential moments by \cite{madaule2011} relying on a work of \cite{maillard:2013}.

 A natural question arising out of the works of \cite{durrett:1983} and \cite{madaule2011} is the following: \emph{where do the point processes based on the scaled displacements converge in the regularly varying case?} The main aim of this article is to show the convergence of this point process sequence and also explicitly identify the limit as a Cox cluster process. We establish that the prediction of \cite{brunet:derrida:2011} on this limit remains true for branching random walk with regularly varying step size even though the finiteness of exponential moments fails to hold. In order to overcome this obstacle, we use a twofold truncation technique based on multivariate extreme value theory.

 We also discuss the superposability properties of our limiting point process in parallel to the recent works of \cite{maillard:2013} and \cite{subag:zeitouni:2014} and confirm the validity of a related prediction of \cite{brunet:derrida:2011} in our setup. As a consequence of our main result, we give explicit formulae for the asymptotic distributions of the properly scaled order and gap statistics from which various problems mentioned in \cite{brunet:derrida:2011} can be investigated. In particular, we recover Theorem~1 of \cite{durrett:1983} in our framework.

This paper is organized as follows. Section~\ref{sec:results} contains the statements of the main result (Theorem~\ref{mainthm}) and its consequences (Theorems~\ref{ssdpppthm1} and \ref{consq1}). Since the proof of Theorem~\ref{mainthm} is long and notationally complicated, we first give a detailed outline of the main steps based on four lemmas in Section~\ref{sec:outline}. These lemmas, and Theorems~\ref{ssdpppthm1} and \ref{consq1} are finally proved in Section~\ref{sec:proofs}.

\section{\bf The Results} \label{sec:results}

We consider point processes as a random elements in the space $\mathscr{M}$ of all Radon point measures on a locally compact and separable metric space $\E$. Here $\mathscr{M}$ is endowed with the vague convergence (denoted by ``$\vconv$''), which is metrizable by the metric
$
\rho(\mu, \nu) = \sum_{i=1}^\infty {2^{-i}} \min \lfb |\mu(h_i)-\nu(h_i)| , 1 \rfb,
$
where $\lsb h_i \rsb_{i \geq 1}$ is a suitably chosen subset (consisting only of Lipschitz functions) of the collection $C_c^+(\E)$ of all non-negative continuous real-valued functions on $\E$ with compact support. $(\mathscr{M}, \rho)$ is a complete and separable metric space. Therefore the standard theory of weak convergence is readily available for point processes and can be characterized by the pointwise convergence of corresponding Laplace functionals on $C_c^+(\E)$ (see Proposition 3.19 in \cite{resnick:1987}). For further details on point processes, see \cite{kallenberg:1986}, \cite{resnick:1987}, \cite{embrechts:kluppelberg:mikosch:1997} and \cite{resnick:2007}.

Because of \eqref{nmr3} and \eqref{nmr4}, we can choose scaling constants $b_n$ such that (see, e.g., \cite{resnick:1987}, \cite{davis:resnick:1985}, \cite{davis:hsing:1995})
\begin{equation}
\mu^n \prob (b_n^{-1} \Xu \in \cdot ) \vconv \nu_\alpha \label{nmr5}
\end{equation}
on $[-\infty,\infty] \setminus \{0\}$, where
\begin{equation} \label{nmr6}
\nu_\alpha(dx) = \alpha p x^{-\alpha -1} \mathbbm{1}_{(0, \infty)}(x) dx + \alpha q (-x)^{-\alpha -1} \mathbbm{1}_{(-\infty,0)}(x) dx.
\end{equation}
Note that one can write $b_n = \mu^{n/\alpha} L_0(\mu^n)$ for some slowly varying function $L_0$. In this paper, conditioned on the survival of the tree, we investigate the asymptotic behaviour of the sequence of point processes defined by
\begin{equation} \label{nmr7}
N_n= \sum_{|\uv|=n} \delta_{b_n^{-1} S_\uv}, \;\; n \geq 1,
\end{equation}
where $S_\uv$ is as in \eqref{nmr2}.

Let $(\Omega,\mathcal F, \prob)$ be the probability space where all the random variables are defined and let $\pstar$ denote the probability obtained by conditioning $\prob$ on the non-extinction of the underlying Galton-Watson tree.  We shall denote by $\exptn$ and $\estar$, the expectation operators with respect to $\prob$ and $\pstar$, respectively. We introduce two sequences of random variables $\{T_l\}_{l \ge 1}$ and $\{j_l\}_{l \ge 1}$ as follows. Suppose $\{T_l\}_{l \ge 1}$ is a sequence of independent and identically distributed positive integer valued random variables with probability mass function
\begin{align}
\gamma(y) := \prob (T_1=y)= \frac{1}{r} \sum_{i=0}^\infty \frac{1}{\mu^i} \prob(Z_i=y), \;\; y \in \mathbb{N}, \label{mr1}
\end{align}
where $r= \sum_{i=0}^\infty \frac{1}{\mu^i} \prob (Z_i>0)$. Let $\{j_l\}_{l \ge 1}$ be a sequence of random variables such that $\summationl \delta_{j_l} \sim PRM(\nu_\alpha)$, where $PRM(\nu_\alpha)$ denotes Poisson random measure with mean measure as $\nu_\alpha$. We also assume that the sequences $\{T_l\}_{l \ge 1}$ and $\{j_l\}_{l \ge 1}$ are independent of each other and are both independent of the martingale limit $W$.

Our main result says that the limiting point process is a Cox cluster process in which a typical Cox point $(rW)^{1/\alpha}j_l$ appears with random multiplicity $T_l$. The clusters appear here due to the strong dependence structure of the displacement random variables $\lsb \Sv : |\uv|=n \rsb$. The randomness in the intensity measure arises from the martingale limit $W$ in contrast to the light tailed case, where similar randomness comes from the derivative martingale limit (see, e.g., Theorem~1.1 in \cite{aidekon2013}). Note also that a $W$-mixture was already present in Theorem 1 of \cite{durrett:1983}; see Remark~\ref{remark:maxima} below.

\begin{thm} \label{mainthm}
With the  assumptions \eqref{nmr1}, \eqref{nmr3}, \eqref{nmr4} and $\{b_n\}$ as in \eqref{nmr5}, under $\pstar$, the sequence of point processes defined in \eqref{nmr7} converges weakly in the space $\mathscr{M}$ of all Radon point measures on $[-\infty, \infty]\setminus \{0\}$ to a Cox cluster process with representation
\begin{equation}
N_* \eql \summationl T_l \delta_{(rW)^{1/\alpha} j_l} \label{nmr13}
\end{equation}
and Laplace functional given by
\begin{align}
& \Psi_{N_*}(g)= \estar \lfb e^{-N_*(g)} \rfb  \nonumber \\
&= \estar \Bigg[ \exp \lsb -W \int_{|x|>0} \sum_{i=0}^\infty \frac{1}{\mu^i}  \exptn ( 1- e^{- Z_i g(x)}) \nu_\alpha{(dx)} \rsb \Bigg] \label{nmr14}
\end{align}
for all $g \in C_c^+([-\infty,\infty]\setminus\{0\})$.
\end{thm}

\subsection{Scale-decorated Poisson point processes} \label{subsec:DPP}

For any point process $\PP$ and any $a >0$, we denote by $\mathbf{s}_a \PP$ the point process obtained by multiplying the atoms of $\PP$ by $a$. The following is an analogue of Definition~1 in \cite{subag:zeitouni:2014} suitable for our framework.

\begin{defn}
\textnormal{A point process $N$ is called a \emph{randomly scaled scale-decorated Poisson point process} (SScDPPP) with intensity measure $\nu$, scale-decoration $\PP$ and random scale $\Theta$ if $N \eql \mathbf{s}_{\Theta} \summation \mathbf{s}_{\lambda_i} \PP_i$, where $\Lambda = \summation \delta_{\lambda_i} $ $\sim PRM(\nu)$ on the space $(0,\infty)$ and $\PP_i$, $i \ge 1$ are independent copies of the point process $\PP$ and are independent of $\Lambda$, and $\Theta$ is a positive random variable independent of $\Lambda$ and $\{\PP_i\}_{i \ge 1}$. We shall denote this by $N \sim SScDPPP (\nu, \PP, \Theta)$. If $\Theta \equiv 1$, we call $N$ a \emph{scale-decorated Poisson point process} (ScDPPP) and denote it by $N \sim ScDPPP(\nu, \PP)$.}
\end{defn}

Our next result establishes that the limiting point process \eqref{nmr13} admits an SScDPPP representation and confirms that a prediction of \cite{brunet:derrida:2011} remains valid in our setup. Moreover, the scale-Laplace functional of $N_*$ (i.e., the left hand side of \eqref{ssdppp1} below) can always be expressed as a multiplicative convolution of an $\alpha\mbox{-Fr\'echet}$ distribution with some measure. This is a scale-analogue of a property investigated in \cite{subag:zeitouni:2014} (see property (SUS) therein).

\begin{thm} \label{ssdpppthm1}
Under the assumptions of Theorem \ref{mainthm}, the limiting point process $N_*  \sim SScDPPP\big( \nu^+_\alpha, T \delta_\varepsilon, (rW)^{1/\alpha}\big)$, where $\varepsilon$ is a $\pm 1$-valued random variable with $\prob(\varepsilon =1)=p$, $\nu^+_\alpha(dx)=\alpha x^{-\alpha -1} dx$ is a measure on $(0,\infty)$ and $T$ is a positive integer valued random variable (independent of $\varepsilon$) with probability mass function \eqref{mr1}. Furthermore, for all $g \in C_c^+([-\infty,\infty]\setminus\{0\})$, $N_*$ satisfies
\begin{align} \label{ssdppp1}
\estar \bigg( \exp \Big\{ - \int g(x/y) N_*(dx) \Big\} \bigg) = \estar \bigg( \Phi_\alpha \Big(c_g\, y\, W^{-1/\alpha} \Big) \bigg), \;\; y>0,
\end{align}
where $\Phi_\alpha$ denotes the distribution function of an $\alpha\mbox{-Fr\'echet}$ random variable, i.e. $\Phi_\alpha(x)= \exp \{ -x^{-\alpha}\}$, $x>0$, and $c_g$ is a positive constant that depends on $g$ but not on $y$.
\end{thm}

The scale-decoration in the SScDPPP representation of $N_*$ is the point process consisting of $T$ many repetitions of the random point $\varepsilon$. This is due to the fact that very few (more precisely, a $W$-mixture of Poisson many) edge random variables survive the scaling by $b_n$ and the surviving ones come with random cluster-sizes that are independent copies of $T$. The presence of $\varepsilon$ in the scale-decoration can be justified by the fact that the surviving edge random variables are positive and negative with probabilities $p$ and $q$, respectively (see \eqref{nmr4}).

\begin{remark}[Superposability]
\textnormal{Let $N^{(i)}_*  = \summationl T^{(i)}_l \delta_{(rW_i)^{1/\alpha} j^{(i)}_l}$, $i=1,2$ be two independent copies of \eqref{nmr13}. Then using Laplace functionals, it can easily be verified that for two positive constants $a_1$ and $a_2$, $\mathbf{s}_{a_1} N^{(1)}_* + \mathbf{s}_{a_2} N^{(2)}_* \sim SScDPPP\big( \nu^+_\alpha, T \delta_\varepsilon, (r (a_1^{\alpha} W_1 + a_2^{\alpha} W_2)^{1/\alpha})\big)$. In particular, when the underlying Galton-Watson tree is a $d$-regular tree (i.e., $Z_1 \equiv d \geq 2$), then the limiting point process is the Poisson cluster process $N_{*,d} \sim ScDPPP\big( \nu^+_\alpha, d^G\delta_\varepsilon\big)$, where $G$ follows a $\mbox{Geometric}(1/d)$ distribution (independently of $\varepsilon$) with probability mass function $\prob (G=k)=(1-1/d)^k d^{-1}$, $k \geq 0$, and $N_{*,d}$ satisfies the superposability property described as follows. If $N_{*,d}^{(i)}$, $i=1,2$ are two independent copies of $N_{*,d}$, then
$
\mathbf{s}_{a_1} N_{*,d}^{(1)} + \mathbf{s}_{a_2} N_{*,d}^{(2)} \eql N_{*,d}
$
for any two positive constants $a_1$, $a_2$ such that $a_1^\alpha + a_2^\alpha =1$.  Following \cite{davydov:molchanov:zuyev:2008}, $N_{*,d}$ can be viewed as an  strictly $\alpha$-stable point process and hence is expected to have an ScDPPP representation (see Section 3 of aforementioned reference). More generally, $\mathbf{s}_{W^{-1/\alpha}}N_*$ is an $\alpha$-stable point process. For similar statements in case of exp-1-stable point processes, see \cite{brunet:derrida:2011}, \cite{maillard:2013} and \cite{subag:zeitouni:2014}.}
\end{remark}

\subsection{{Order and gap statistics}}
The point process convergence in Theorem~\ref{mainthm} helps us to derive some properties of the order and gap statistics; see also \cite{ramola:majumder:schehr:2014} for related works on branching Brownian motion. Let $M_n^{(k)}$ denote the $k^{th}$ upper order statistic coming from the $n^{th}$ generation, $G_n^{(k)}= M_n^{(k)} - M_n^{(k+1)}$ be the $k^{th}$ gap statistic and $M'_n:= \min_{|\uv|=n} \Sv$ be the minima. In order to study the asymptotic properties of these statistics, we need a few more notations as described below. We denote by $\pi$ a partition of an integer $l$ of the form $l=i_1 y_{1} + i_2 y_{2} + \cdots + i_{|\pi|} y_{|\pi|}$, where each $i_j$ repeats $y_j$ many times in the partition, and  $i_1< i_2 <\cdots < i_{|\pi|}$. Here $|\cdot|$ denotes the number of distinct elements in a partition. Let $\Pi_l$ be the set of all such partitions of the integer $l$.

\begin{thm}\label{consq1}
With the assumptions of Theorem \ref{mainthm} and $\{b_n\}$ as in \eqref{nmr5}, the following asymptotic properties hold.
\begin{itemize}

\item[\textnormal{(a)} ] \textnormal{(Minima)} For all $x >0$, $$\displaystyle{\lim_{n \to \infty} \pstar \lfb M'_n >- b_n x \rfb = \estar \lfb \exp \lsb -r W q x^{- \alpha} \rsb \rfb}.$$

\item[\textnormal{(b)} ] \textnormal{($k^{th}$ upper order statistic)} For all $x >0$,
\begin{align}
& \lim_{n \to \infty} \pstar \lfb M_n^{(k)}  \le b_n x \rfb = \estar \Bigg( \exp \bigg\{- rW p x^{-\alpha} \bigg\} \Bigg) \label{con2p2}\\
& + \sum_{l=1}^{k-1} \sum_{\pi \in \Pi_l} \estar \Bigg[ \prod_{j=1}^{|\pi|}  \bigg( \Big( r W p  x^{-\alpha} \gamma(i_j) \Big)^{y_{j}} \exp \left\{ -r W p x^{-\alpha} \gamma(i_j)\right\}  \bigg) \Bigg].  \nonumber
\end{align}

\item[\textnormal{(c)} ] \textnormal{(Joint distribution of $k^{th}$ and $(k+1)^{th}$ upper order statistics)} For all $(u,v)$ such that $ 0 < u < v $,
\begin{align} \label{jointorder}
& \lim_{n \to \infty} \pstar ( M_n^{(k+1)} \le b_n u, M_n^{(k)} \le b_n v) = \estar \Big( \xi_{0, (u, \infty]}(W) \Big) \\
& +  \sum_{j=1}^{k} \estar \Big( \xi_{0, (v, \infty]} (W) \xi_{j, (u,v]}(W)\Big)  + \sum_{l=1}^{k-1} \sum_{j=0}^{k-l} \estar \Big( \xi_{l, (v, \infty]} (W) \xi_{j, (u,v]}(W)\Big), \nonumber
\end{align}
where for all $l \ge 0$ and for all $ A \subset [-\infty, \infty] \setminus \{0\}$ such that $\nu_\alpha(A) < \infty$,
\begin{equation*}
\xi_{l,A} (W) := \begin{cases}  e^{  -r W \nu_\alpha(A) }  \h \h \hspace{2.1cm}   \mbox{ if } l=0, \\
  \displaystyle{\sum_{\pi \in \Pi_l} \prod_{j=1}^{|\pi|} (rW \nu_\alpha(A) \gamma(i_j))^{y_j}  \frac{1}{y_j !} e^{- r W \nu_\alpha(A) \gamma(i_j) }} \hspace{.6cm} \mbox{ if }  l \ge 1. \end{cases}
\end{equation*}

\item[\textnormal{(d)} ] \textnormal{($k^{th}$ gap statistic)} Let $L : \mathbb{R}^+\times \mathbb{R}^+\to \mathbb{R}^+$ be the map $L(u,v)= v-u$. Then $\pstar ( b_n^{-1} G_n^{(k)} \in \cdot) \to \zeta_k\circ L^{-1}$ where $\zeta_k$ is a probability measure on $\mathbb{R}^+ \times \mathbb{R}^+$ with joint cumulative distribution function \eqref{jointorder}.
\end{itemize}
The second term in \eqref{con2p2} and the third term in \eqref{jointorder} are both interpreted as zero when $k=1$.
\end{thm}

\begin{remark}[Maxima] \label{remark:maxima}
\textnormal{Note that putting $k=1$ in \eqref{con2p2}, we recover Theorem~1 of \cite{durrett:1983} in our framework. By Theorem 8.2 (Page 15) of \cite{harris:1963} (see also \cite{athreya:ney:1972}, Theorem 2, Page 29), the limiting distribution function of the scaled maxima of the $n^{th}$ generation can be written as
\begin{equation}
 \lim_{n \to \infty} \pstar \lfb M_n^{(1)}  \le b_n x \rfb = \phi(rpx^{-\alpha}), \;\; x >0,\label{limit:maxima:alt:form}
\end{equation}
where $\phi$ is the unique (up to a scale-change) completely monotone function on $\mathbb{R}^{+}$ satisfying
\begin{equation}
 \phi(z)=f(\phi(z/\mu))  \label{eqn:functional}
\end{equation}
with $f$ being the probability generating function of the branching random variable $Z_1$.}
\end{remark}

\begin{example}[Maxima for geometric branching] \textnormal{Suppose that the offspring distribution of the underlying branching process is geometric with parameter $b \in (0,1)$ and probability mass function $\prob (Z_1=k)= b (1-b)^{k-1}$, $k \ge 1$. It is easy to check that the completely monotone function $\phi(u)= \frac{1}{1+du}$, $u >0$ satisfies the functional equation \eqref{eqn:functional} for any scaling constant $d>0$. Therefore using \eqref{limit:maxima:alt:form} and the fact that $E(W)=1$ (a consequence of Kesten-Stigum condition~\eqref{nmr1}; see \cite{kesten:stigum:1966a}), it follows that $d=1$ and
\begin{align*}
\lim_{ n \to \infty} \prob ( M_n \le b_n x) = \frac{1 - b}{1 - b + px^{-\alpha}}, \;\; x>0.
\end{align*}}
\end{example}

\section{\bf Outline of Proof of Theorem \ref{mainthm}} \label{sec:outline}

In this section, we outline the main result's proof, which is based on a twofold truncation technique using extreme value theory. We attain this via four lemmas, whose proofs will be given in the next section. For ease of presentation, we shall use Ulam-Harris labeling system described recursively as follows. The $i^{th}$ descendant of the root $o$ is denoted by $i$ and $j^{th}$ descendant of an $(n-1)^{th}$ generation vertex $(i_1, \ldots, i_{n-1})$ is denoted by $(i_1, \ldots , i_{n-1},j)$. We abuse the notation and denote an edge joining an $(n-1)^{th}$ generation vertex and an $n^{th}$ generation vertex using the same label as the latter vertex. Such an edge is assumed to belong to the $n^{th}$ generation. Let $D_n$ denote the vertices (and hence edges because of the abuse of notation) in the $n^{th}$ generation and $C_n = \cup_{i=1}^{n} D_i$ denote the vertices (as well as edges) up to the $n^{th}$ generation of the underlying Galton-Watson tree. With these notations, we describe below the mains steps of the proof of Theorem \ref{mainthm}.

\subsection{One large jump} \label{onelargejump}

Following \cite{durrett:1983}, it is easy to see that with very high probability, for every vertex $\uv \in D_n$, among all $\uu \in \Iv$ at most one edge random variable $\Xu$ will be large enough    to survive the scaling by $b_n$. Hence we can expect that the asymptotic behavior of $N_n$ will be same as that of
\begin{equation} \label{pmr1}
\tilde{N}_n= \sum_{|\uv|=n} \sum_{\uu \in \Iv} \delta_{b_n^{-1} \Xu}.
\end{equation}
More precisely, we shall establish the following lemma.

\begin{lemma} \label{lemma1}
Under the assumptions of Theorem~\ref{mainthm}, for every $\epsilon >0$,
\begin{equation} \label{pmr2}
\limsup_{n \to \infty} \pstar \lfb \rho \lfb N_n, \tilde{N}_n \rfb > \epsilon \rfb =0,
\end{equation}
where $\rho$ is the vague metric introduced at the beginning of Section~\ref{sec:results}.
\end{lemma}
This lemma formalizes the well-known principle of one large jump (see, e.g., Steps 3 and 4 in Section~2 of \cite{durrett:1983}) at the level of point processes and it can be shown by molding the proof of Theorem 3.1 in \cite{resnick:samorodnitsky:2004}. Because of Lemma~\ref{lemma1}, it is enough to investigate the weak convergence of \eqref{pmr1}, which is much easier compared to that of \eqref{nmr7}.

\subsection{Cutting the tree}\label{cutting}

The first truncation is a standard one that has been used in branching random walks. First fix a positive integer $K$. Taking $n>K$, look at the tree $\T$ up to the $n^{th}$ generation and cut it at $(n-K)^{th}$ generation keeping last $K$ generations alive; see Figure 1. This means that after cutting the tree, we will be left with a forest containing $K$ generations of $|D_{n-K}|$ many independent (under $\prob$) Galton-Watson trees with roots being the vertices at the $(n-K)^{th}$ generation of the original tree $\T$ and the same offspring distribution as before. We label the new sub-trees in this forest as $\lsb \T_j \rsb_{j=1}^{|D_{n-K}|}$.
\begin{center}

\begin{tikzpicture}

 \draw (0,0) node [anchor = south]{$o$};

 \draw[ dotted,->] (0,0) --  (-2,-2);
 \draw [ dotted, ->] (0,0) --  (2,-2);


 \draw[dotted, ->] (-2,-2) --  (-5,-4) ;
 \draw[dotted, ->] (-2,-2) --  (-2.5,-4);
 \draw[dotted,  ->] (-2,-2) --  (-0,-4);

 \draw[dotted,  ->] (2,-2) --  (2.5,-4);
 \draw[dotted, ->] (2,-2) --  (5,-4);

 \draw[dotted]  (-6,-4) -- (6,-4);

 \draw  (-5,-4) node[left]{$\T_1$};
 \draw  (-2.5,-4) node[left]{$\T_2$};
 \draw (-0,-4) node[left]{$\T_3$};
 \draw  (2.5,-4) node[left]{$\T_4$};
 \draw  (5,-4) node[left]{$\T_5$};


 \draw[ ->] (-5,-4) -- (-5.5,-6);
 \draw[ ->] (-5,-4) -- (-4.5,-6);

 \draw[ ->] (-2.5,-4) -- (-3,-6);
 \draw[ ->] (-2.5,-4) --  (-2.5,-6)  ;
 \draw[ ->] (-2.5,-4) -- (-2,-6);

 \draw[ ->] (2.5,-4) -- (2.5,-6);

 \draw[ ->] (5,-4) -- (4.5,-6);
 \draw[ ->] (5,-4) -- (5,-6);
 \draw[ ->] (5,-4) -- (5.5,-6);

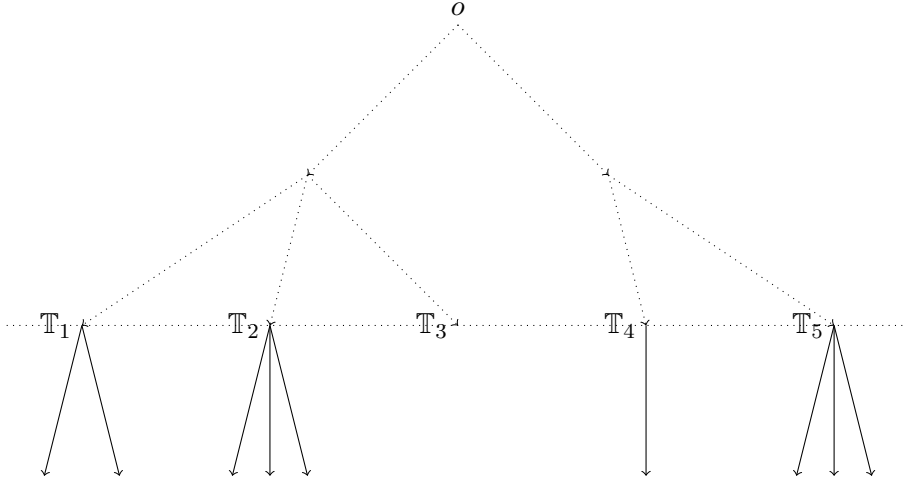
 \draw (0,-6.5) node[anchor=north]{ \small{Figure 1  : Cutting the Galton-Watson tree ($n=3, K=1$) at generation $2$.}};

 \end{tikzpicture}

\end{center}

Each vertex $\uv$ in the $n^{th}$ generation of the original tree $\T$ belongs to the $K^{th}$ generation of some sub-tree $\T_j$ and we denote by $\Iv^K$ the unique geodesic path from the root of $\T_j$ to the vertex $\uv$. We introduce another point process generated by the i.i.d.\ heavy-tailed random variables attached to the edges of the forest as follows:
\begin{equation} \label{pmr3}
\tnk := \sum_{|\uv|=n} \sum_{\uu \in \Iv^K} \delta_{b_n^{-1}\Xu},
\end{equation}
where $|\uv|$ denotes the generation of $|\uv|$ in the original tree $\T$. The following lemma asserts that as long as $K$ is large, \eqref{pmr3} is a good approximation of \eqref{pmr1}.

\begin{lemma} \label{lemma2}
Under the assumptions of Theorem~\ref{mainthm}, for every $\epsilon > 0$,
\begin{equation} \label{pmr4}
\lim_{K \to \infty} \limsup_{n \to \infty} \pstar \lfb \rho \lfb \tilde{N}_n, \tnk \rfb > \epsilon \rfb =0.
\end{equation}
\end{lemma}

In light of the above lemma, it is enough to find the weak limit of \eqref{pmr3} as $n \to \infty$ keeping $K$ fixed, and then letting $K \to \infty$. This can be achieved with the help of another truncation as mentioned below.

\subsection{Pruning the forest}\label{pruning}

This is the second truncation step, which is also quite standard in branching process theory. Fix an integer $K>0$ and for each edge $\uu$ in the forest $\cup_{j=1}^{|D_{n-K}|} \T_j$, define $A_\uu$ to be the number of descendants of $\uu$ at $n^{th}$ generation of $\T$. Fix another integer $B>1$ large enough so that $\mu_B := \exptn (Z_1^{(B)}) > 1$, where $Z_1^{(B)} := Z_1 \mathbbm{1}(Z_1 \le B) + B \mathbbm{1}(Z_1 > B)$. We modify the forest according to the pruning algorithm mentioned below (see also Figure 2).

\begin{enumerate}

\item[P1.] Start with the sub-tree $\T_1$ and look at its root.

\item[P2.] If the root has more than $B$ many children (edges), then keep the first $B$ many edges according to our labeling, and delete the others and their descendants. If the number of children (edges) of the root is less than or equal to $B$, then do nothing.

\item[P3.] Now we can have at most $B$ many vertices in the first generation of the sub-tree $\T_1$. Repeat Step P2 for children (edges) of each of these vertices. Continue with this algorithm up to the children (edges) of the $(K-1)^{th}$ generation vertices (of the sub-tree $\T_1$).

\item[P4.] Repeat Steps P2 and P3 for the other sub-trees $\T_2, \ldots, \T_{|D_{n-K}|}$.

\end{enumerate}

\begin{center}

\begin{tikzpicture}

\draw[dotted]  (-6,-4) -- (6,-4);

 \draw  (-5,-4) node[left]{$\T_1$};
 \draw  (-2.5,-4) node[left]{$\T_2$};
 \draw (0,-4) node[left]{$\T_3$};
 \draw  (2.5,-4) node[left]{$\T_4$};
 \draw  (5,-4) node[left]{$\T_5$};


 \draw[ ->] (-5,-4) -- (-5.5,-6);
 \draw[ ->] (-5,-4) -- (-4.5,-6);

 \draw[ ->] (-2.5,-4) -- (-3,-6);
 \draw[ ->] (-2.5,-4) --  (-2.5,-6)  ;

 \draw[ ->] (2.5,-4) -- (2.5,-6);

 \draw[ ->] (5,-4) -- (4.5,-6);
 \draw[ ->] (5,-4) -- (5,-6);

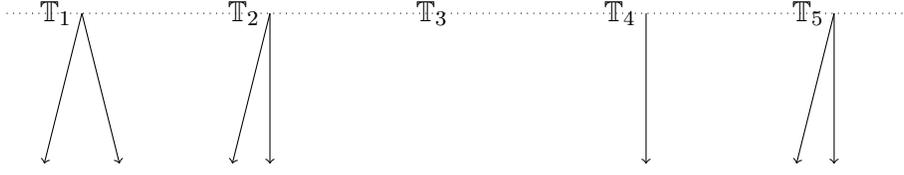
 \draw (0,-6.5) node[anchor=north]{ \small{Figure 2: Pruning of the forest obtained in Figure 1 with $B=2$.}};

 \end{tikzpicture}

\end{center}

Note that under $\prob$, these $|D_{n-K}|$ many pruned sub-trees are independent copies of a Galton-Watson tree (up to the $K^{th}$ generation) with a bounded branching random variable $Z_1^{(B)}$. For each $j$, we denote by $\T^{(B)}_j$ the pruned version of $\T_j$. For each edge $\uu$ in $\cup_{j=1}^{|D_{n-K}|} \T^{(B)}_j$, we define $\aub$ to be the number of descendants of $\uu$ in the $K^{th}$ generation of the corresponding pruned sub-tree. Observe that for every vertex $\uu$ at the $i^{th}$ generation of any sub-tree $\T^{(B)}_j$, $\aub$ is equal is distribution to $Z_{K-i}^{(B)}$, where $\lsb Z_i^{(B)} \rsb_{i \ge 0}$ denotes a branching process with $Z_0^{(B)} \equiv 1$ and branching random variable $Z_1^{(B)}$. For each $i = 1, 2. \ldots, K$, we denote by $D_{n-K+i}^{(B)}$ the union of all $i^{th}$ generation vertices (as well as edges) from the pruned sub-trees $\T^{(B)}_j$, $j=1,2,\ldots, |D_{n-K}|$. We introduce another point process as follows.
\begin{equation} \label{pmr7}
\nkb := \sum_{\uv \in D_n^{(B)}} \sum_{\uu \in \Iv} \delta_{b_n^{-1} \Xu}.
\end{equation}

The point processes $\tnk$ and $\nkb$ are not simple point processes since both of them have alternative representations as given below.
\begin{equation} \label{pmr8}
\tnk = \sum_{i=0}^{K-1} \sum_{\uu \in D_{n-i}} A_\uu \delta_{b_n^{-1} \Xu}
\end{equation}
and
\begin{equation} \label{pmr9}
\nkb = \sum_{i=0}^{K-1} \sum_{\uu \in D_{n-i}^{(B)}} \aub \delta_{b_n^{-1} \Xu}.
\end{equation}

The set of all trees upto $K^{th}$ generation becomes a finite set due to pruning. This helps in the computation of the limit \eqref{eq:tree_wald} below. The next lemma encompasses this second truncation step and reduces our work to computation of weak limit of \eqref{pmr9} obtained by letting $n \to \infty$, and then $B \to \infty$, and finally $K \to \infty$.

\begin{lemma} \label{lemma3}
  Under the assumptions of Theorem~\ref{mainthm}, for each fixed positive integer $K>1$ and for all $\epsilon > 0$,
\begin{equation} \label{pmr10}
\lim_{B \to \infty} \limsup_{n \to \infty} \pstar \lfb \rho \lfb \tnk, \nkb \rfb >\epsilon \rfb=0.
\end{equation}
\end{lemma}

\subsection{Computation of Weak limit}
We shall compute the weak limit of  \eqref{pmr9} by investigating its Laplace functional. This is carried out using  a conditioning argument and extreme value theory. The conditioning argument helps us to exploit the independence in the underlying tree structure. This results in the following lemma, which should be regarded as the key step in proving Theorem~\ref{mainthm}.

\begin{lemma} \label{lemma4}
Under the conditions of Theorem~\ref{mainthm}, the following weak convergence results hold in the space $\mathscr{M}$ of all Radon point measures on $[-\infty, \infty]\setminus \{0\}$ under the measure $\pstar$.
\begin{itemize}
\item[(a)] For each positive integer $K$ and each integer $B>1$ with $\mu_B >1$, there exists a point process $N_*^{(K,B)}$ such that $\nkb \Rightarrow N_*^{(K,B)}$ as $n \to \infty$.
\item[(b)] For each positive integer $K$, there exists a point process $N_*^{(K)}$ such that $N_*^{(K,B)} \Rightarrow N_*^{(K)}$ as $B \to \infty$.
\item[(c)] As $K \to \infty$, $N_*^{(K)} \Rightarrow N_*$.
\end{itemize}
\end{lemma}

For detailed descriptions of the point processes $N_*^{(K,B)}$ and $N_*^{(K)}$, see Section~\ref{sec:proofs} below.

\subsection{Proof of Theorem \ref{mainthm}} \label{prooflimitpp}

It is easy to check that \eqref{nmr13} is $\pstar$-almost surely Radon and hence is a random element of $\mathscr{M}$ with $\mathbb{E}=[-\infty,\infty]\setminus\{0\}$. To compute the Laplace functional of $N_*$, take any $g \in C_c^+([-\infty, \infty]\setminus \{0\})$ and observe that
\begin{equation}
\estar (e^{ - N_*(g) }) = \estar \Big[ \estar \big(  e^ {  - N_*(g)}\big| W\big) \Big] = \estar \Big[ \estar \big(  e^ {  - N(f)}\big| W\big) \Big] , \label{eq:alt_lap_fl_Nstar}
\end{equation}
where $f$ is the function $f(t,x)=t g(x)$ defined on $\mathbb{N} \times ([-\infty, \infty]\setminus \{0\})$ and $N$ is the Cox process
$
N = \summationl \delta_{(T_l, \,(rW)^{1/\alpha}j_l)}.
$
Using Propositions~3.6 and 3.8 in \cite{resnick:1987}, we get
$$
\estar  \big( e^{- N(f)} \big| W \big) = \exp \Big \{ -rW   \int_{|x|>0}   \exptn (1- e^{-f(T_1,x)})\nu_\alpha(dx) \Big \},
$$
which can be shown to be equal to the random quantity inside the expectation in \eqref{nmr14}. Therefore, the second part of Theorem~\ref{mainthm} follows from \eqref{eq:alt_lap_fl_Nstar}.

Using Lemmas~\ref{lemma2}, \ref{lemma3} and \ref{lemma4} and applying twice a standard converging together argument (see, e.g., Theorem~3.5 in \cite{resnick:2007}), it follows that under $\pstar$,
\begin{equation}
\tilde{N}_n \Rightarrow N_* \mbox{   as $n \to \infty$}, \label{conv_of_Ntilde}
\end{equation}
from which the weak convergence in Theorem~\ref{mainthm} follows by a simple application of Theorem~3.4 of \cite{resnick:2007} combined with  Lemma~\ref{lemma1} above.

\section{ \bf Rest of the Proofs} \label{sec:proofs}

Throughout this section $\PT$ will denote the probability obtained by conditioning $\pstar$ on the whole Galton-Watson tree $\T$ and $\E$ will denote the space $[-\infty, \infty] \setminus \{0\}$. Also we will use the notation $\surv$ to denote the event that the Galton-Watson tree survives.

\subsection{ Proof of Lemma \ref{lemma2}}

Let $g \in C_c^+(\E)$ with $\mbox{support}(g) \subseteq \lsb x : |x|>\delta \rsb$, for some $\delta>0$. By definition of the vague convergence, it is enough to show that for all   $\epsilon>0$,
\begin{equation} \label{lm3p1}
\lim_{K \to \infty} \limsup_{n \to \infty} \pstar \lfb | \tilde{N}_n(g)- \tilde{N}_n^K(g)| > \epsilon \rfb =0.
\end{equation}
Recall that $C_{n-K} = \{\uv : |\uv| \le n-K\}$. Define $B_{n,K}$ to be the event that all the random variables in the collection $\{X_\uu: \uu \in C_{n-K}\}$ are less than $b_n \delta /2$ in modulus .
We claim that $\lim_{K \to \infty} \limsup_{ n \to \infty} \pstar \lfb B_{n,K}^c \rfb=0$, which will follow provided we show that
\begin{equation} \label{lm3p5}
\lim_{K \to \infty} \limsup_{n \to \infty} \PT \lfb B_{n,K}^c \rfb =0 \mbox{   for $\pstar$-almost all $\T$.}
\end{equation}

To this end, note that conditioned on the tree $\T$, $\sum_{\uu \in C_{n-K}} \delta_{b_n^{-1} |\Xu|} (\theta, \infty]$ follows a $\mbox{Binomial}(|C_{n-K}|, \prob (|\Xu|> b_n \theta))$ distribution, and for each $K \geq 1$, the following $\pstar$-almost sure convergence holds:
\begin{eqnarray*}
|C_{n-K}| \prob \lfb |\Xu| > b_n \theta \rfb &=& \frac{|C_{n-K}|}{\mu^{n-K}} \frac{1}{\mu^K} \mu^n \prob \lfb |\Xu| > b_n \theta \rfb \nonumber \\
& \stackrel{\scriptscriptstyle{n \to \infty}}{\to} & \frac{\mu}{\mu-1} W \frac{1}{\mu^K} \theta^{- \alpha} =: \lambda(\theta, K).
\end{eqnarray*}
Therefore, for all $K \geq 1$, $\sum_{\uu \in C_{n-K}} \delta_{b_n^{-1} |\Xu|} (\theta, \infty] \Rightarrow \mathscr{P} \sim \mbox{Poisson}(\lambda(\theta,K))$ as $n \to \infty$ under $\PT$ for $\pstar$-almost all $\T$. Because of Kesten-Stigum condition \eqref{nmr1}, $\lambda(\theta,K) \to 0$ $\pstar$-almost surely as $K \to \infty$ . In particular, we get $\lim_{n \to \infty} \PT \lfb B_{n,K}^c\rfb =\PT (\mathscr{P} >1)$, which tends to zero as $K \to \infty$ for $\pstar$-almost all $\T$ and hence~\eqref{lm3p5} holds. To finish the proof from here, observe that $\PT \lfb |\tilde{N}_n(g) - \tilde{N}_n^K(g)| > \epsilon, B_{n,K} \rfb \equiv 0$ since the support of $g$ is contained in $\lsb x : |x| > \delta \rsb$. Hence \eqref{lm3p1} follows immediately from~\eqref{lm3p5}.

\subsection{Proof of Lemma \ref{lemma3}}
Let $g\in C_c^+(\E)$ be as in proof of Lemma~\ref{lemma2} with $\mbox{support}(g) \subseteq \{x: |x|>\delta\}$ and $\|g\|_\infty:=\sup_{x \in \mathbb{E}} |g(x)| < \infty$. To show~\eqref{pmr10}, it is enough to show that for such $g\in C_c^+(\E)$ and $\epsilon>0$,

\begin{equation*}
\lim_{B \to \infty} \limsup_{n \to \infty} \pstar \bigg[  \Big| \tilde{N}_n^K(g)- \tilde{N}_n^{(K,B)}(g)  \Big| > \epsilon \bigg] =0.
\end{equation*}
Noting that the points from the point process $\tilde{N}_n^{(K,B)}$ are contained in the point process $\tilde{N}_n^K$, and using \eqref{pmr8} and \eqref{pmr9}, we have $ | \tnk (g) - \tnkb (g) |= \sum_{i=0}^{K-1} \left(S^{(1)}_{i,n,B} + S^{(2)}_{i,n,B} \right)$, where $S^{(1)}_{i,n,B}=\sum_{\uu \in D_{n-i}} (A_\uu-\aub) g(b_n^{-1} \Xu) $ and $S^{(2)}_{i,n,B}=\sum_{\uu \in D_{n-i} \setminus \dnib} \aub g(b_n^{-1} \Xu) $. Since $\pstar(\cdot) \leq (\prob(\surv))^{-1}\prob(\cdot)$, it is enough to show that for each $i$, both $S^{(1)}_{i,n,B}$ and $S^{(2)}_{i,n,B}$ are negligible under $\prob$.

To this end, fix $0 \leq i \leq K-1$ and $\eta>0$. Using Markov's inequality, Wald's identity and the bound $|g| \leq \|g\|_\infty \bm_{[-\infty,-\delta)\cup(\delta,\infty]}$, we get
\begin{align*}
\prob\left(S^{(1)}_{i,n,B} > \eta\right) &\leq \frac{1}{\eta} \exptn (Z_{n-i}) \exptn\left(A_\uu-\aub\right) \exptn \left(g(b_n^{-1} X_\uu)\right)\\
                                                                & \leq \frac{1}{\eta \mu^i} \|g\|_\infty\exptn\left(A_\uu-\aub\right)\mu^n\prob\left(|X_\uu|>b_n \delta\right),
\end{align*}
from which first letting $n \to \infty$ based on \eqref{nmr5} and then letting $B \to \infty$, it follows that $\lim_{B \to \infty} \limsup_{n \to \infty} \prob\left(S^{(1)}_{i,n,B} > \eta\right)=0$. We can deal with $S^{(2)}_{i,n,B}$ in a similar fashion and obtain
\begin{align*}
\prob\left(S^{(2)}_{i,n,B} > \eta\right) &\leq \frac{\|g\|_\infty}{\eta}\exptn \left(\aub\right)\prob (|X_\uu|> b_n \delta) \exptn ( |D_{n-i}| - |\dnib|)\\
                                                                &\leq \frac{\|g\|_\infty}{\eta}\exptn (A_\uu)\mu^n \prob (|X|> b_n \delta) \left(\frac{1}{\mu^i} - \frac{\mu_B^K}{\mu_B^i} \mu^{-K}\right).
\end{align*}
Therefore, $\lim_{B \to \infty} \limsup_{n \to \infty} \prob\left(S^{(2)}_{i,n,B} > \eta\right)=0$. This suffices.

\subsection{Proof of Lemma \ref{lemma4}}

We shall first establish (a). To this end, we introduce an event $S_{n-K}$ which is empty when $|\dnk| =0$, and on $(|\dnk|>0)$, it is the event that there is at least one infinite tree rooted at the $(n-K)^{th}$ generation of the underlying Galton-Watson tree. Using $\1_{\surv} = \1_{(|\dnk|>0)} \1_{S_{n-K}}$
and $d \pstar = \prob(\surv)^{-1} \1_\surv d \prob$, we get that
for every $g \in C_c^+([-\infty, \infty] \setminus \{0\})$,
\begin{align}
&\estar \bigg[ \exp \{- \tnkb(g)\} \bigg]
= \frac{1}{\prob(\mathcal{S})} \exptn \bigg[\1_{\mathcal{S}} \exp \{- \tnkb(g)\} \bigg] \nonumber \\
& =  \frac{1}{\prob(\mathcal{S})} \exptn \bigg[\1_{(|\dnk|>0)}  \exp \{- \tnkb(g)\} \bigg] \label{eq:any_label} \\
& \h  -  \frac{1}{\prob(\mathcal{S})} \exptn \bigg[\1_{(|\dnk|>0)} \1_{S_{n-K}^c} \exp \{- \tnkb(g)\} \bigg] \nonumber
\end{align}
It is easy to see that the second term of the above is bounded by
\begin{align*}
\frac{1}{\prob (\surv)} \exptn (\bm_{\{Z_{n-K} >0 \} } \bm_{S_{n-K}^c}) = \frac{1}{\prob(\surv)} \exptn \Big( \bm_{\{Z_{n-K} >0\}} p_e^{Z_{n-K}} \Big),
\end{align*}
where $p_e = \prob(\surv^c)$ denotes the probability of extinction of the underlying Galton-Watson tree. Since our tree is supercritical, $p_e <1$ and $\prob$-almost surely, $Z_{n-K} \to \infty$ as $n \to \infty$ on the event $\mathcal{S}$. Using the $\prob$-almost sure convergence $\bm_{\{Z_{n-K} >0\}} \to \bm_\surv $, it is easy to see that second term of \eqref{eq:any_label} tends to zero as $n \to \infty$.

 It is enough to study asymptotics of~\eqref{eq:any_label} . We shall first condition on the $\mathcal{F}_{n-K}$ and then use the fact that conditioned on $\mathcal{F}_{n-K}$, $\{(A_\uu, \Xu)_{\uu \in \T_j} : j=1 , \ldots, |D_{n-K}|\}$ is a collection of independent and identically distributed random variables under $\prob$. This will imply that \eqref{eq:any_label} is equal to
\aln{
& \frac{1}{\prob(\mathcal{S})} \exptn \bigg[ \exptn \bigg( \1_{(|\dnk|>0)}  \exp \bigg\{- \sum_{j=1}^{|\dnk|} \sum_{\uu \in \T_j^{(B)}} A_{\uu}^{(B)} g(b_n^\inv \Xu) \bigg\} \bigg| \mathcal{F}_{n-K} \bigg) \bigg] \nonumber \\
&=\frac{1}{\prob(\mathcal{S})} \exptn \Bigg[ \1_{(|\dnk|>0)} \Bigg(\exptn \bigg(   \exp \bigg\{-  \sum_{\uu \in \tob} A_{\uu}^{(B)} g(b_n^\inv \Xu) \bigg\}  \bigg) \Bigg)^{|\dnk|} \Bigg]. \label{eq:alt_proof1}
}

Define $\mathcal{T}_{K,B}$ to be the set of all rooted trees having at most $K$ generations and each vertex having at most $B$ branches. Then the inner expectation in \eqref{eq:alt_proof1} can be written as
\alns{
& 1 - \frac{1}{\mu^n} \sum_{\bft \in \caltkb} \prob(\tob = \bft)  \exptn \bigg[ \mu^n \bigg( 1- e^{ - \sum_{\uu \in \bft} A_\uu^{(B)} g(b_n^\inv \Xu)  } \bigg)  \bigg| \tob = \bft \bigg].
}
For each fixed $ \bft \in \mathcal{T}_{K,B}$, we shall first compute the limit (as $n \to \infty$) of the conditional expectation in the above expression. To this end, we shall denote the edges of the tree by the pair $(i,j)$, which will indicate the $j^{th}$ edge in the $i^{th}$ generation of $\bft$. The corresponding $A^{(B)}_\uu$ and $X_\uu$ will be denoted by $A^{(B)}_{(i,j)}$ and $X_{(i,j)}$, respectively. Also $|\bft|_i$ will denote the total number of vertices in the $i^{th}$ generation of $\bft$ and $|\bft|:=\sum_{i=1}^K |t|_i$. With these notations, we have
\alns{
&\exptn \bigg[ \mu^n \bigg( 1- e^{ - \sum_{\uu \in \bft} A_\uu^{(B)} g(b_n^\inv \Xu)  } \bigg)  \bigg| \tob = \bft \bigg]\\
&=  \int \bigg( 1- \exp \Big\{ - \sum_{i =1}^K \sum_{j=1}^{\cardi{\bft}} A_{(i,j)}^{(B)} g( x_{(i,j)})  \Big\} \bigg)  \mu^n  \prob \bigg( b_n^\inv \tilde{X} \in d\tilde{x} \bigg| \tob = \bft \bigg)
}
where, conditioned on the event $\{\tob = \bft\}$, $\tilde{X}$ denotes the vector
$$(X_{(1,1)}, \ldots, X_{(1, |\bft|_1)}, \ldots, X_{(i,1)}, \ldots, X_{(i,\cardi{\bft})}, \ldots, X_{(K,1)}, \ldots, X_{(K, |\bft|_K)})$$
 and
$$\tilde{A}^{(B)}=(A_{(1, 1)}^{(B)}, \ldots, A_{(1, |\bft|_1)}^{(B)}, \ldots A^{(B)}_{(i,1)}, \ldots, A_{(i, |\bft|_i)}^{(B)}  \ldots, A_{(K,1)}^{(B)}, \ldots, A_{(K, |\bft|_K)}^{(B)}).$$

 Using the fact that $\Xu$'s are independent and identically distributed random variables satisfying \eqref{nmr5}, we get
\aln{
\mu^n  \prob(b_n^\inv \tilde{X} \in \cdot | \tob = \bft) \vconv \tau_\bft(\cdot) = \sum_{i=1}^K \sum_{j=1}^{\cardi{\bft}} \tau^\bft_{(i,j)}(\cdot) \label{eq:mult_reg_var}
}
on the space $[-\infty, \infty]^{|\bft|} \setminus \{\mathbf{0}\}$, where for all $1 \le i \le K$ and $1 \le j \le \cardi{\bft}$
\alns{
\tau^\bft_{(i,j)} : = \delta_0 \times \ldots \times \delta_0 \times \underbrace{\nu_\alpha}_{(|\bft_1| + \cdots + |\bft|_{i-1} + j)^{th} \mbox{ position}} \times \delta_0 \times \ldots \times \delta_0 .
}
Note that on the event $\{\tob = \bft\}$, $\tilde{A}^{(B)}$ is a deterministic vector completely specified by $\bft$.  Hence using \eqref{eq:mult_reg_var}, we get that
\alns{
& \int \bigg( 1- \exp \Big\{ - \sum_{i =1}^K \sum_{j=1}^{\cardi{\bft}} A_{(i,j)}^{(B)} g (x_{(i,j)})  \Big\} \bigg) \mu^n  \prob \bigg( b_n^\inv \tilde{X} \in d\tilde{x} \bigg| \tob = \bft \bigg) \\
&\to\sum_{i=1}^K \sum_{j=1}^{\cardi{\bft}} \int_{|x|>0} \bigg( 1- \exp \Big\{ -  A_{(i,j)}^{(B)} g( x)  \Big\} \bigg) \nu_\alpha(d x).
}
Since $\mathcal{T}_{K,B}$ is a finite set, it follows that as $n \to \infty$,
\begin{align}
& \sum_{\bft} \prob(\tob = \bft)  \int \bigg( 1- \exp \Big\{ - \sum_{i =1}^K \sum_{j=1}^{\cardi{\bft}} A_{(i,j)}^{(B)} g( x_{(i,j)})  \Big\} \bigg) \nonumber  \\
& \hspace{5cm}\mu^n  \prob \bigg( b_n^\inv \tilde{X} \in d\tilde{x} \bigg| \tob = \bft \bigg)  \nonumber\\
& \to \sum_{i=1}^K \exptn \bigg[  \sum_{j=1}^{\cardi{\tob}} \int_{|x|>0} \bigg( 1- \exp \Big\{ -  A_{(i,j)}^{(B)} g( x)  \Big\} \bigg) \nu_\alpha(d x) \bigg].\label{eq:tree_wald}
\end{align}

Let $Z_i^{(B)}$ denote the number of particles in the $i^{th}$ generation of the Galton-Watson process with branching random variable $Z_1^{(B)}$. For every fixed $i$, $|\mathbb T_1^{(B)}|_i$ and $\{A_{(i,j)}^{(B)}: j\ge 1\}$ are independent, and $\{A_{(i,j)}^{(B)}: j\ge 1\}$ is a sequence of i.i.d.\ random variables with distribution as that of $Z_{K-i}^{(B)}$. Using Wald's identity we get that~\eqref{eq:tree_wald} equals
\alns{
\sum_{i=1}^K \mu_B^i \int_{|x|>0} \bigg( 1- \exp\bigg\{- Z_{K-i}^{(B)} g(x) \bigg\}\bigg) \nu_\alpha(d x).
}
Combining this with the almost sure convergence $\mu^{-n} |D_{n-K}| \to \mu^{-K} W$ (as $n \to \infty$) it follows that
\alns{
 & \Bigg(\exptn \bigg(   \exp \bigg\{-  \sum_{\uu \in \tob} A_{\uu}^{(B)} g(b_n^\inv \Xu) \bigg\}  \bigg) \Bigg)^{|\dnk|} \\
&=   \Bigg( 1 - \frac{1}{\mu^n} \sum_{\bft \in \caltkb} \prob(\tob = \bft) \\
& \h \exptn \bigg[ \mu^n \bigg( 1- \exp \Big\{ - \sum_{\uu \in \bft} A_\uu^{(B)} g(b_n^\inv \Xu)  \Big\} \bigg)  \bigg| \tob = \bft \bigg] \Bigg)^{\mu^n \frac{|\dnk|}{\mu^n}} \\
& \to  \exp \bigg\{ - \frac{1}{\mu^K} W \sum_{i=1}^K \mu_B^i \int_{|x|>0} \exptn \bigg( 1- \exp\bigg\{- Z_{K-i}^{(B)} g(x) \bigg\}\bigg) \nu_\alpha(d x) \bigg\}
}
almost surely as $n \to \infty$. Note that $\1_{(|\dnk|>0)} \to \1_\surv$ as $n \to \infty$. Therefore an application of dominated convergence theorem yields, \eqref{eq:alt_proof1} converges to
\alns{
& \estar \Bigg[ \exp \Bigg\{ - W \frac{1}{\mu^K}  \sum_{i=1}^K \mu_B^i \exptn \bigg[  \int_{|x|>0} \bigg( 1- e^{ - Z^{(B)}_{K-i} g( x)}   \bigg) \nu_\alpha(d x) \Bigg\} \Bigg]  \\
& =  \estar \Bigg[ \exp \Bigg\{ - W \frac{\mu_B^K}{\mu^K}  \sum_{i=0}^{K-1} \mu_B^{i} \exptn \bigg[  \int_{|x|>0} \bigg( 1- e^{  - Z^{(B)}_{i} g( x)  } \bigg) \nu_\alpha(d x) \Bigg\} \Bigg].
}
This can easily be shown (using an argument similar to the one used in Subsection~\ref{prooflimitpp}) to be the Laplace functional of the point process
$$
N_*^{(K,B)} := \summationl T_l^{(K,B)} \delta_{(r_{K,B} (\frac{\mu_B}{\mu})^K W)^{1/\alpha}j_l} \, ,
$$
where $r_{K,B}= \sum_{i=0}^{K-1} \frac{1}{\mu_B^i} \prob (Z_i^{(B)} > 0)$, $\{j_l\}_{l \ge 1}$ is a sequence of random variables such that $\summationl \delta_{j_l} \sim PRM(\nu_\alpha)$ and $\{T_l^{(K,B)}\}_{l \ge 1}$ is a sequence of i.i.d.\ random variables (independent of $\{j_l\}_{l \geq 1}$ and $W$) with probability mass function
$$
\prob(T_1 =y) = \frac{1}{r_{K,B}} \sum_{i=0}^{K-1} \frac{1}{\mu_B^i} \prob (Z_i^{(B)}=y), \;y \in \mathbb{N}.
$$
Thus (a) follows using Theorem~5.2 in \cite{resnick:2007}.

To establish (b), fix a positive integer $K$ and observe that applying dominated convergence theorem as $B \to \infty$, the Laplace functional of $\nkb$ can be shown to converge to that of
$$
N_*^{(K)} := \summationl T_l^{(K)} \delta_{(r_K W)^{1/\alpha} j_l} \, ,
$$
where $r_{K}= \sum_{i=0}^{K-1} \frac{1}{\mu^i} \prob (Z_i > 0)$ and $\{T_l^{(K)}\}_{l \ge 1}$ is a sequence of i.i.d.\ random variables (independent of $\{j_l\}_{l \geq 1}$ and $W$) with probability mass function
$$
\prob(T_1 =y) = \frac{1}{r_{K}} \sum_{i=0}^{K-1} \frac{1}{\mu^i} \prob (Z_i=y), \;y \in \mathbb{N}.
$$
In a similar fashion, (c) can be shown. This completes the proof of Lemma~\ref{lemma4}.

\subsection{ Proof of Lemma \ref{lemma1}}

To show~\eqref{pmr2}, it is enough to take a Lipschitz function $g\in C_c^+(\E)$ (with Lipschitz constant $\|g\|$ and $\mbox{support}(g) \subseteq \lsb x : |x| > \delta \rsb$ for some $\delta>0$) and show that for every $\epsilon >0$,

\begin{equation} \label{lm4p1}
\lim_{n \to \infty} \pstar \lfb |N_n(g)-\tilde{N}_n(g)| > \epsilon \rfb =0.
\end{equation}
This will be attained by slightly revamping the proof of the convergence in (3.14) of \cite{resnick:samorodnitsky:2004}. Some of the estimates used therein will not work for us mainly because we are dealing with general regularly varying random variables as opposed to stable ones with an inbuilt Poissonian structure. This hurdle will be overcome by use of Potter's bound and a mild modification of the event AMO$(\theta)$ defined in page~201 of the aforementioned reference.

For every $\theta>0$, let $A_n(\theta) $ denote the event that for all $\uv \in D_n$, at most one of the random variables in the collection $\{X_\uu: \uu \in I_\uv\}$ is bigger than $b_n \theta/n$ in absolute value. We claim that
\begin{equation}\label{lm4p4}
\lim_{n \to \infty} \pstar \lfb A_n(\theta)^c \rfb=0.
\end{equation}
 As in the proof of Lemma~\ref{lemma2}, this follows easily if we can establish that
$\lim_{n \to \infty} \PT \lfb A_n(\theta)^c \rfb =0$  for $\pstar$-almost all $\T$.
To this end, observe that conditioned on the tree $\T$, $\sum_{\uu \in \Iv} \delta_{|\Xu|} (b_n \frac{\theta}{n}, \infty]$ $\sim \mbox{Binomial}(n, \prob(n |\Xu|> b_n \theta))$ for each $\uv \in D_n$. Hence using Potter's bound (see, e.g., Proposition 0.8(ii) in \cite{resnick:1987}), \eqref{nmr1}, \eqref{nmr5}, and the fact that $P(U_n \geq 2) = O(n^2 p_n^2)$ for any $U_n \sim \mbox{Binomial}(n, p_n)$ with $p_n=\mathrm{o}(1/n)$, we get

\begin{align*}
\PT \lfb A_n(\theta)^c \rfb  &\le  \sum_{|\uv|=n} \PT \lfb \sum_{\uu \in \Iv} \delta_{|X_\uu|} ( b_n \frac{\theta}{n}, \infty]  \ge 2 \rfb \\
                                                 &\leq (\mbox{const})  |D_n| n^2 \big(\prob \lfb n |\Xu| > b_n \theta \rfb\big)^2 \to 0
\end{align*}
$\pstar$-almost surely as $n \to \infty$.

In light of \eqref{lm4p4}, to prove \eqref{lm4p1}, it is enough to show that
\begin{equation}
\lim_{n \to \infty} \pstar\lfb |N_n(g)- \tilde{N}_n(g)| > \epsilon , A_n(\theta)\rfb = 0.\label{limit_needed_to_prove}
\end{equation}
Let us fix $0 <\theta < \delta/2$. On the event $A_n(\theta)$, define $T_\uv$ to be the the largest (in absolute value) summand  in $\sum_{\uu \in \Iv} \Xu$.  Note that
\begin{align}
&\PT ( |N_n (g) - \tilde{N}_n(g) | > \epsilon, A_n(\theta)) \nonumber\\
&\le \PT \bigg( \sum_{|\uv| = n} | g(b_n^\inv \Sv) - \sum_{\uu \in \Iv} g(b_n^\inv \Xu)| > \epsilon, A_n(\theta) \bigg)  \nonumber\\
&= \PT \bigg( \sum_{|\uv| = n} | g(b_n^\inv \Sv) -  g(b_n^\inv T_\uv)| > \epsilon, A_n(\theta) \bigg) .
\label{eq:one large_max}
\end{align}
For every $\uv \in D_n$, $|\Sv - T_\uv| < b_n \theta < b_n \delta/2$ on the event $A_n(\theta)$.
So for fixed $\uv \in D_n$, $|g(b_n^\inv \Sv) - g(b_n^\inv T_\uv)| >0$ only if $b_n^\inv T_\uv >\delta/2$.
Using the fact that $g$ is Lipschitz, we get that for every $\uv$, $|g(b_n^\inv \Sv) - g(b_n^\inv T_\uv)| \le \|g\| b_n^\inv |\Sv - T_\uv| \le \|g\| \theta$. So ~\eqref{eq:one large_max} can be bounded by
\alns{
& \PT \bigg(   \|g\| \theta  \sum_{|\uv| =n} \delta_{b_n^\inv T_\uv} (\delta/2, \infty] >\epsilon  \bigg) \\
& \le \PT \bigg(   \|g\| \theta  \sum_{|\uv| =n} \sum_{\uu \in  \Iv}\delta_{b_n^\inv \Xu}\Big( [-\infty, \delta/2) \cup (\delta/2, \infty] \Big)>\epsilon  \bigg) \\
& = \PT \bigg( \|g\| \theta  \tilde{N}_n\Big( [-\infty, \delta/2) \cup (\delta/2, \infty] \Big)>\epsilon \bigg).
}
Unconditioning the above expression and using the fact that $\tilde{N}_n \Rightarrow N_*$ as $n \to \infty$, we get that left hand side of ~\eqref{limit_needed_to_prove} converges to
\alns{
\pstar \bigg( \|g\| \theta   N_*\Big( [-\infty, \delta/2) \cup (\delta/2, \infty] \Big)>\epsilon \bigg).
}
 Now let $\theta \to 0$ to get \eqref{lm4p1}.

%

\begin{remark}
\textnormal{The proof of Lemma~\ref{lemma1} uses \eqref{conv_of_Ntilde}, which is a consequence of Lemmas~\ref{lemma2}, \ref{lemma3} and \ref{lemma4}; see Subsection~\ref{prooflimitpp}. However, this is not a problem because the latter lemmas are proved without using Lemma~\ref{lemma1}.}
\end{remark}

\subsection{{ Proof of Theorem \ref{ssdpppthm1}}}
Take two independent sequences of random variables $\{\varepsilon_i\}_{i \geq 1}$ and $\{\lambda_i\}_{i \geq 1}$ such that $\sum_{i=1}^\infty \delta_{\lambda_i} \sim PRM(\nu^+_\alpha)$ and $\varepsilon_1, \varepsilon_2, \ldots $ are i.i.d.\ random variables with same distribution as that of $\varepsilon$. Straightforward applications of Propositions~5.2 and 5.3 of \cite{resnick:2007} yield that $\sum_{i=1}^\infty \delta_{\varepsilon_i\lambda_i} \sim PRM(\nu_\alpha)$, which, together with \eqref{nmr13}, gives the SScDPPP representation of $N_*$.

To show the second part of this theorem, we follow the computation of Laplace functional in Subsection \ref{prooflimitpp} along with the scaling property of $\nu_\alpha$ and express the left hand side of \eqref{ssdppp1} as
\begin{align*}
& \estar \Bigg( \exp \bigg\{ -W\int_{|x|>0} \sum_{i=0}^\infty \frac{1}{\mu^i} \exptn \Big( 1- e^{-Z_i f(y^{-1} x)} \Big) \nu_\alpha(dx) \bigg\} \Bigg) \\
&=\estar \Bigg( \exp \bigg\{ - y^{-\alpha}W\int_{|x|>0} \sum_{i=0}^\infty \frac{1}{\mu^i}  \exptn \Big( 1- e^{-Z_i g( x)} \Big) \nu_\alpha(dx) \bigg\} \Bigg),
\end{align*}
which equals the right hand side with
\begin{align*}
c_g = \left(\int_{|x|>0} \sum_{i=0}^\infty \frac{1}{\mu^i} \exptn \Big( 1- e^{-Z_i g( x)} \Big) \nu_\alpha(dx)\right)^{-1/\alpha} >0.
\end{align*}

\subsection{Proof of Theorem~\ref{consq1}}

Using Theorem~3.1 and Theorem~3.2 of \cite{resnick:2007} and Theorem~\ref{mainthm} above, it transpires that $N_n([ -\infty,-x])$ converges weakly to $N_*([-\infty,-x])$ under $\pstar$. Therefore, for each $x>0$,
$$
\pstar\big(M'_n > - b_n x\big) = \pstar\big(N_n([-\infty,-x])=0\big) \to \pstar\big(N_*([-\infty,-x])=0\big),
$$
from which (a) follows because $T_l > 0$ for all $l \ge 1$ and this implies
\begin{align*}
\pstar\big(N_*([-\infty,-x])=0 |W\big) &= \pstar\bigg(\summationl \delta_{(rW)^{1/\alpha}j_l}([-\infty,-x])=0\Big|W\bigg)\\
                                                                &=  \exp \lsb -rWq x^{-\alpha} \rsb.
\end{align*}

The $k=1$ case of (b) follows similarly from the weak convergence of $N_n((x,\infty])$ to $N_*((x,\infty])$ under $\pstar$. For $k \geq 2$, using the same weak convergence, we get
\begin{align}
& \lim_{n\to \infty} \pstar \lfb M_n^{(k)}  \le b_n x\rfb =\lim_{n\to \infty} \pstar ( N_n((x, \infty]) \le k-1) \nonumber \\
&= \estar \left(\exp \lsb -rWp x^{-\alpha} \rsb\right) + \sum_{l=1}^{k-1} \pstar ( N_*((x, \infty)) = l).\label{con2p3}
\end{align}
We need to show that the second term of \eqref{con2p3} is same as that of \eqref{con2p2}. To this end, considering the marked point process
$N = \summationl \delta_{(T_l, (rW)^{1/\alpha}j_l)} \sim PRM\big(rW(\gamma \otimes \nu_\alpha)\big)$ conditioned on $W$, and analyzing exactly how each event $(N_*((x, \infty)) = l)$ can occur, the second term in \eqref{con2p3} becomes
\begin{align*}
& \sum_{l=1}^{k-1} \sum_{ \pi \in \Pi_l} \pstar \lfb \bigcap_{j=1}^{|\pi|}  \lsb N( \lsb i_j \rsb \times (x, \infty] ) =y_{j} \rsb \rfb  \\
&= \sum_{l=1}^{k-1} \sum_{\pi \in \Pi_l} \estar \Bigg[ \prod_{j=1}^{|\pi|}  \bigg( \Big( r W p  x^{-\alpha} \gamma(i_j) \Big)^{y_{j}}\frac{1}{y_{j} !}  \exp \Big\{ -r W p x^{-\alpha} \gamma(i_j) \Big\} \bigg) \Bigg].
\end{align*}
This establishes (b).

In order to verify (c), we need a similar (but slightly tedious) calculation as in the proof of (b) based on the following observation: for $0 < u < v$,
\begin{align*}
&\pstar \big( M_n^{(k+1)} \le b_n u, M_n^{(k)} \le b_n v\big)\\
&\; = \pstar\big(N_n((u,\infty])=0\big)  + \pstar\big(N_n((v,\infty])=0, 1 \leq N_n((u,v]) \leq k\big)\\
&\;\;\;\;\;+ \pstar\big(1 \leq N_n((v,\infty]) \leq k-1, N_n((u,\infty]) \leq k\big).
\end{align*}
Finally, (d) follows from (c) using continuous mapping theorem (see, e.g., Theorem~3.1 in \cite{resnick:2007}).

\noindent \section*{Acknowledgements} The authors would like to thank  Gregory Schehr for asking a question that resulted in Theorem~\ref{consq1}~(d) above. The authors are also very grateful to the anonymous referees for their suggestions which have improved the paper significantly.

\end{document}